\documentclass[12pt, reqno]{amsart}

\usepackage{amsmath, amssymb, amsthm}

\DeclareMathOperator{\PG}{PG}

\DeclareMathOperator{\trace}{Tr}

\DeclareMathOperator{\rad}{Rad}

\numberwithin{equation}{section}

\newcommand{\F}{{\mathbb F}}
\newcommand{\cC}{{\mathcal C}}
\newcommand{\cF}{{\mathcal F}}
\newcommand{\cK}{{\mathcal K}}
\newcommand{\floor}[1]{\left\lfloor #1\right\rfloor}

\newcommand{\abs}[1]{\left\vert #1\right\vert}
\newcommand{\ontop}[2]{\genfrac{}{}{0pt}{}{#1}{#2}}
\newcommand{\ba}{{\mathbf a}}

\newcommand{\bu}{{\mathbf u}}
\newcommand{\bbu}{{\mathbf{\bar u}}}
\newcommand{\bv}{{\mathbf v}}
\newcommand{\by}{{\mathbf y}}
\newcommand{\bw}{{\mathbf w}}

\newcommand{\bz}{{\mathbf z}}
\newcommand{\bzero}{{\mathbf 0}}
\newcommand{\bone}{{\mathbf 1}}
\newcommand{\doesdivide}{{\,|\,}}

\newtheorem{theorem}{Theorem}[section]
\newtheorem*{theorem*}{Theorem}
\newtheorem{lemma}[theorem]{Lemma}

\newtheorem{example}[theorem]{Example}

\begin{document}
\title[Maximal Arcs in Desarguesian Planes] {On Mathon's Construction of
  Maximal Arcs in Desarguesian Planes. II.}
\author{Frank Fiedler, Ka Hin Leung, Qing Xiang}
\address{Department of Mathematical Sciences, University of Delaware,
  Newark, DE 19716, USA, email: fiedler@math.udel.edu}
\address{Department of Mathematics, National University of Singapore, Kent
  Ridge, Singapore 119260, email: matlkh@nus.edu.sg}
\address{Department of Mathematical Sciences, University of Delaware,
  Newark, DE 19716, USA, email: xiang@math.udel.edu}


\begin{abstract}
 In a recent paper \cite{Mat02}, Mathon gives a new construction
 of maximal arcs which generalizes the construction of Denniston.
 In relation to this construction, Mathon asks the question of determining the
 largest degree of a non-Denniston maximal arc arising from his
 new construction. In this paper, we give a nearly complete answer to
 this problem. Specifically, we prove that when $m\geq 5$ and $m\neq 9$,
 the largest $d$ of a non-Denniston maximal arc of degree $2^d$ in $\PG(2,2^m)$
 generated by a $\{p,1\}$-map is $(\floor{\frac{m}{2}}+1)$. This confirms
 our conjecture in \cite{FLX}. For $\{p,q\}$-maps, we prove that if $m\geq 7$ and $m\neq 9$,
 then the largest $d$ of a non-Denniston maximal arc of degree $2^d$ in $\PG(2,2^m)$
 generated by a $\{p,q\}$-map is either $\floor{\frac{m}{2}}+1$ or $\floor{\frac{m}{2}}+2$.
\end{abstract}

\maketitle

\section{Introduction}
\label{section:introduction} Let $\PG(2,q)$ denote the
desarguesian projective plane of order $q$, where $q$ is a prime
power, and let $k\geq 1$, $n\geq 2$ be integers. A set $\cK$ of
$k$ points in $\PG(2,q)$ is called a $(k,n)$-\emph{arc} if no
$n+1$ points of $\cK$ are collinear.  The integer $n$ is called
the \emph{degree} of the arc $\cK$. Let $P$ be a point of a
$(k,n)$-arc $\cK$. Each of the $q+1$ lines through $P$ contains at
most $n-1$ points of $\cK$. Therefore
\[
k\le 1+(q+1)(n-1).
\]
The $(k,n)$-arc $\cK$ is said to be \emph{maximal} if $k$ attains
this upper bound, that is, $k=q (n-1)+n$. In this case, every line
of $\PG(2,q)$ that contains a point of $\cK$ has to intersect it
in exactly $n$ points. Therefore the degree $n$ of a maximal arc
$\cK$ in $\PG(2,q)$ must divide $q$.

In the case where $q=2^m$, maximal arcs of degree 2 in $\PG(2,q)$
are usually called \emph{hyperovals}. A classical example of a
hyperoval in $\PG(2,2^m)$ is a non-degenerate conic (i.e.,
non-singular quadric in $\PG(2,2^m)$) plus its nucleus. There is
an extensive literature devoted to ovals and hyperovals, see a
recent survey in \cite{Pen}. The study of maximal arcs of degree
greater than two was started by Barlotti \cite{B} in 1955. At the
beginning, maximal arcs were studied as extremal objects in finite
geometry and coding theory. Later it was discovered that maximal
arcs can give rise to many interesting incidence structures such
as partial geometries, resolvable Steiner 2-designs
(\cite{Thas74}, \cite{wallis}). The constructions of Thas
\cite{Thas74, Thas80} also show connections between maximal arcs
with ovoids, quadrics and polar spaces. Of course maximal arcs can
also give rise to two-weight codes and strongly regular graphs
since they are two-intersection sets in $\PG(2,q)$. For these
reasons maximal arcs occupy a very special place in finite
geometry, design theory and coding theory.

For $q=2^m$, Denniston \cite{Den69} constructed maximal arcs of
degree $2^d$ in $\PG(2,2^m)$ for every $d$, $1\leq d\leq m$. Thas
\cite{Thas74}, \cite{Thas80} also gave two other constructions of
maximal arcs in $\PG(2,2^m)$ of certain degrees when $m$ is even.
For odd prime power $q$, Ball, Blokhuis and Mazzocca \cite{bbm}
proved that maximal arcs of degree $n$ do not exist in $\PG(2,q)$,
when $n<q$. Recently Mathon \cite{Mat02} presented a new
construction of maximal arcs which generalizes the construction of
Denniston. In the following, we will briefly describe the
constructions of Denniston and Mathon of maximal arcs.

Let $Q(x,y)=a x^2 + h x y + b y^2$ be an irreducible quadratic
form over $\F_{2^m}$ (that is, $\trace(\frac{a b}{h^2})=1$, where
$\trace$ is the trace from $\F_{2^m}$ to $\F_2$). Let
$A$ be an additive subgroup of $\F_{2^m}$ and let
$(x,y,z)$ be right-normalized homogeneous coordinates in
$\PG(2,2^m)$. Then
\begin{equation}
\label{definition:DennistonArc}
\cK=\{ (x,y,1)\in\PG(2,2^m)\mid Q(x,y)\in A\}
\end{equation}
is a maximal arc of degree $\abs{A}$. This is Denniston's construction of
maximal arcs \cite{Den69}. We may decompose $\cK$ as
\[
\cK=\cup_{\lambda\in A} F_{\lambda},
\]
where for each $\lambda\in A\setminus\{0\}$, $F_{\lambda}=\{(x,y,1)\mid
Q(x,y)=\lambda\}$ is a non-degenerate conic, and $F_{0}=\{(0,0,1)\}$
contains one point only. Note that the point $(0,0,1)$ is the common nucleus
of the conics $F_{\lambda}$, $\lambda\in A\setminus\{0\}$. The arc $\cK$ in
(\ref{definition:DennistonArc}), and those projectively equivalent to $\cK$
are called \emph{Denniston maximal arcs}.

Now let $\cC$ be the set of conics
\[
 F_{\alpha,\beta,\lambda} = \{ (x,y,z)\in\PG(2,2^m)\mid
 \alpha x^2 + x y + \beta y^2 + \lambda z^2 = 0\},
\]
where $\lambda\in\F_{2^m}\cup\{\infty\}$ and
$\alpha,\beta\in\F_{2^m}^*$ such that $\alpha x^2 + x + \beta$ is
irreducible over $\F_{2^m}$. Note that $F_0:=F_{\alpha, \beta, 0}=
\{(0,0,1)\}$ is the common nucleus of the non-degenerate conics in
$\cC$, and $F_\infty:=F_{\alpha, \beta, \infty}$ is the line at
infinity $z=0$. Given two non-degenerate conics
$F_{\alpha,\beta,\lambda}$ and $F_{\alpha',\beta',\lambda'}$ in
$\cC$ with $\lambda\neq\lambda'$, Mathon \cite{Mat02} defined a
composition
\begin{equation}
\label{definition:composition} F_{\alpha,\beta,\lambda}\oplus
F_{\alpha',\beta',\lambda'}=F_{\alpha\oplus\alpha',\beta\oplus\beta',\lambda\oplus\lambda'},
\end{equation}
where
\[
\alpha\oplus\alpha' = \frac{\alpha\lambda+\alpha'\lambda'}{\lambda+\lambda'},
\quad
\beta\oplus\beta' = \frac{\beta\lambda+\beta'\lambda'}{\lambda+\lambda'},
\quad
\lambda\oplus \lambda'=\lambda+\lambda'.
\]
A subset of non-degenerate conics of $\cC$ that is closed under the above
composition is called a {\it closed} set of conics, and such a set must
contain $2^d-1$ conics for some $d$, $1\le d\le m$
(\cite[Corollary~2.3]{Mat02}). Mathon \cite{Mat02} showed that closed sets
of conics can be used to construct maximal arcs.

\begin{theorem}[\protect{\cite[Theorem~2.4]{Mat02}}]
\label{theorem:Mathon:closed}
Let $\cF\subset \cC$ be a closed set of $2^d-1$ non-degenerate conics with a
common nucleus $F_0$ in $\PG(2,2^m)$, $1\le d\le m$. Then the set of points
of all conics in $\cF$ together with $F_0$ form a maximal
$(2^{m+d}-2^m+2^d,2^d)$-arc $\cK$ in $\PG(2,2^m)$.
\end{theorem}

The construction in the Theorem~\ref{theorem:Mathon:closed} clearly contains
Denniston's construction of maximal arcs as a special case. Let $A$ be an
additive subgroup of $\F_{2^m}$, let $a,b,h\in \F_{2^m}$ be fixed such that
$\trace(\frac{a b}{h^2})=1$, and let $\cF=\{F_{ah^{-1},bh^{-1},\lambda
  h^{-1}}\in \cC\mid \lambda\in A\setminus\{0\}\}$. Then $\cF$ is clearly
closed with respect to the composition in (\ref{definition:composition}),
and the maximal arc obtained via Theorem~\ref{theorem:Mathon:closed} from
$\cF$ is exactly the Denniston arc in (\ref{definition:DennistonArc}).

Let $\cF\subset\cC$ be a closed set of $(2^d-1)$ non-degenerate
conics, and let
\[
A^*=\{\lambda\mid \text{there exist } \alpha, \beta\in
\F_{2^m}^*\; \text{ such that } F_{\alpha,\beta,\lambda}\in
\cF\}.\] Then $A:=A^*\cup\{0\}$ is an additive subgroup of
$\F_{2^m}$. Moreover, for each $\lambda\in A^*$ there corresponds
a unique conic $F_{\alpha,\beta,\lambda}$ in $\cF$ (otherwise,
$F_0\in \cF$, a contradiction), hence $\alpha$ and $\beta$ in the
indices of $F_{\alpha,\beta,\lambda}$ can be interpreted as
functional values of some functions $p:A\to \F_{2^m}$ and $q:A\to
\F_{2^m}$, respectively. Since $\cF$ is closed under the
composition defined in (\ref{definition:composition}), we have
\begin{align*}
  p(\lambda+\lambda')(\lambda+\lambda') & =
  p(\lambda)\lambda+p(\lambda')\lambda',\\
  q(\lambda+\lambda')(\lambda+\lambda') & =
  q(\lambda)\lambda+q(\lambda')\lambda'
\end{align*}
for all $\lambda,\lambda'\in A$. Thus, the maps $\bar
p:A\to\F_{2^m}$ and $\bar q:A\to\F_{2^m}$ defined respectively by
$\bar p(\lambda)= p(\lambda)\lambda$ and $\bar q(\lambda)=
q(\lambda)\lambda$ are linear on $A$. Since $A$ is an
$\F_2$-subspace of $\F_{2^m}$, we can extend $\bar p$ and $\bar q$
linearly to $\F_{2^m}$, and we denote the extended maps still by
$\bar p$ and $\bar q$. Now that $\bar p$ and $\bar q$ are both
linear on $\F_{2^m}$, there exist linearized polynomials
$\sum_{i=0}^{m-1}c_ix^{2^i}$ and $\sum_{i=0}^{m-1}d_ix^{2^i}$ in
$\F_{2^m}[x]$ such that for all $a\in\F_{2^m}$, ${\bar
p}(a)=\sum_{i=0}^{m-1}c_ia^{2^i}$ and ${\bar
  q}(a)=\sum_{i=0}^{m-1}d_ia^{2^i}$. Furthermore, by ``division algorithm''
(c.f. \cite[Proposition 3.1]{FLX}), there exist linearized polynomials
$L(x)=\sum_{i=0}^{d-1}a_ix^{2^i}$ and $M(x)=\sum_{i=0}^{d-1}b_ix^{2^i}$ in
$\F_{2^m}[x]$ such that ${\bar p}(\lambda)=L(\lambda)$ and ${\bar
  q}(\lambda)=M(\lambda)$ for all $\lambda\in A$. This shows that each
closed set $\cF\subset\cC$ of $(2^d-1)$ conics can be written in the form
\[
\{F_{\frac{L(\lambda)}{\lambda}, \frac{M(\lambda)}{\lambda}, \lambda}\mid
\lambda\in A\setminus\{0\}\},
\]
where $A$ is some additive subgroup of $\F_{2^m}$ of size $2^d$, and $L(x),
M(x)\in \F_{2^m}[x]$ are given above.

\begin{theorem}[\protect{\cite[Theorem~2.5]{Mat02}}]
  \label{theorem:Mathon:pq}
  Let $p(x)=\sum_{i=0}^{d-1} a_i x^{2^i-1}\in\F_{2^m}[x]$ and
  $q(x)=\sum_{i=0}^{d-1} b_i x^{2^i-1}\in\F_{2^m}[x]$ be polynomials with
  coefficients in $\F_{2^m}$. For an additive subgroup $A$ of order $2^d$ of
  $\F_{2^m}$ let $\cF=\{F_{p(\lambda),q(\lambda),\lambda}\mid\lambda\in A\setminus\{0\}\}$
  be a set of conics with a common nucleus $F_0$.  If
  $\trace(p(\lambda)q(\lambda))=1$ for every $\lambda\in A\setminus\{0\}$,
  then $\cF$ is a closed subset of $\cC$  and the set of points on all conics in $\cF$
  together with $F_0$ forms a maximal $(2^{m+d}-2^m+2^d,2^d)$-arc $\cK$ in $\PG(2,2^m)$. If both
  $p(x)$, $q(x)$ have $d\le 2$, then $\cK$ is a Denniston
  maximal arc.
\end{theorem}

We will call maximal arcs generated by polynomials as in the above
theorem \emph{maximal arcs generated by $\{p,q\}$-maps}. Mathon
posed several problems related to the construction in
Theorem~\ref{theorem:Mathon:pq} at the end of his paper
\cite{Mat02}. The third problem he posed is: What is the largest
$d$ of a non-Denniston maximal arc of degree $2^d$ in $\PG(2,2^m)$
generated by a $\{p,q\}$-map via Theorem~\ref{theorem:Mathon:pq}?
We give a nearly complete answer to this problem in this paper
(see details below). The techniques we use are algebraic.
Polynomials over finite fields play an important role throughout
our investigation. Combinatorial and linear algebraic tools are
used to study these polynomials in this paper. We hope that these
techniques will find more applications in finite geometry and
combinatorial designs.

Our main results are summarized as follows. In
Section~\ref{section:p-1}, we prove that if $m\geq 5$ and $m\neq
9$, then the largest degree of a non-Denniston maximal arc in
$\PG(2,2^m)$ generated by a $\{p,1\}$-map is less than or equal to
$2^{\floor{\frac{m}{2}}+1}$. On the other hand, known
constructions in \cite{Mat02}, \cite{HamMat}, \cite{FLX} show that
there are always $\{p,1\}$-maps that generate non-Denniston
maximal arcs in $\PG(2,2^m)$ of degree $2^{\floor{\frac{m}{2}}+1}$
when $m\geq 5$. Therefore, for $\{p,1\}$-maps, we have a complete
answer to Mathon's question mentioned above. That is, when $m\geq
5$ and $m\neq 9$, the largest $d$ of a non-Denniston maximal arc
of degree $2^d$ in $\PG(2,2^m)$ generated by a $\{p,1\}$-map via
Theorem~\ref{theorem:Mathon:pq} is $\floor{\frac{m}{2}}+1$. This
confirms our conjecture in \cite{FLX}. In
Section~\ref{section:p-q} we try to extend this result to
$\{p,q\}$-maps. We prove that if $m\ge 7$ and $m\neq 9$, then the
largest degree of a non-Denniston maximal arc in $\PG(2,2^m)$
generated by a $\{p,q\}$-map is less than or equal to
$2^{\floor{\frac{m}{2}}+2}$. However, at present we are not able
to find a construction of $\{p,q\}$-maps to produce (via
Theorem\ref{theorem:Mathon:pq}) a non-Denniston maximal arc in
$\PG(2,2^m)$ of degree $2^{\floor{\frac{m}{2}}+2}$. Therefore our
upper bound together with previously known constructions in
\cite{Mat02}, \cite{HamMat}, \cite{FLX}, yields that for $m\geq 7$
and $m\neq 9$, the largest $d$ of a non-Denniston maximal arc of
degree $2^d$ in $\PG(2,2^m)$ generated by a $\{p,q\}$-map is
either $\floor{\frac{m}{2}}+1$ or $\floor{\frac{m}{2}}+2$.

\section{The Largest Degree of non-Denniston Maximal Arcs generated by $\{p,1\}$-Maps}
\label{section:p-1}

We first prove the following theorem, which establishes the upper
bound mentioned in Section~\ref{section:introduction} on the
largest degree of non-Denniston maximal arcs generated by a
$\{p,1\}$-map.

\begin{theorem}
  \label{theorem:p1}
  Let $A$ be an additive subgroup of size $2^d$ in $\F_{2^m}$, and let
  $p(x)=\sum_{i=0}^{d-1} a_ix^{2^i-1}\in\F_{2^m}[x]$. Assume that $m\ge
  5$ but $m\neq 9$, and $m>d>\frac{m}{2}+1$. If $\trace(p(\lambda))=1$
  for all $\lambda\in A\setminus\{ 0\}$, then $a_{2}=a_{3}=\cdots
  =a_{d-1}=0$. That is, $p(x)$ is linear and the maximal arc obtained
  from the $\{p,1\}$-map via Theorem~\ref{theorem:Mathon:pq} is a
  Denniston maximal arc.
\end{theorem}

In order to prove this theorem we need some preparation. For
convenience, let $r=m-d$. We will represent the $\F_2$-subspace
$A$ of $\F_{2^m}$ as the intersection of $r$ hyperplanes, say
\[A=\{ x\in\F_{2^m}\mid \trace(\mu_i x)=0, 1\le i\le r\},\]
where $\mu_i\in\F_{2^m}^*$ are linearly independent over $\F_{2}$.  Thus,
the defining equation for $A$ is
\[
\prod_{i=1}^r (1+\trace(\mu_i x))=1.
\]
The key to the proof of Theorem~\ref{theorem:p1} is to study the
polynomial $ \prod_{i=1}^r (1+\trace(\mu_i x))$, where
$\trace(\mu_ix)=\sum_{j=0}^{m-1}\mu_i^{2^j}x^{2^j}$ is a
polynomial in $\F_{2^m}[x]$. We define $S(x)$ to be the polynomial
of degree less than or equal to $2^m-1$ such that
\[ S(x) \equiv \prod_{i=1}^r (1+\trace(\mu_i x)) \pmod {x^{2^m}-x}.\]
For $s\ge 1$ and $m-1\ge i_1>i_2>\cdots >i_s\ge 0$, we use $c(i_1,
i_2, \ldots, i_s)$ to denote the coefficient of
$x^{2^{i_1}+2^{i_2}+\cdots +2^{i_s}}$ in $S(x)$.  It is clear that
$c(i_1, i_2, \ldots, i_s)$ is zero if $ s> r$. Moreover, as
$S(x)^2 \equiv S(x) \pmod {x^{2^m}-x}$, we see that when $s\leq
r$,
\begin{equation}
  \label{equation:coefficient:shift}
  c(i_1,i_2,\ldots , i_s)^2  = \left\{\begin{array}{ll}
c(i_1+1,i_2+1\ldots, i_s+1) & \mbox{ if } i_1<m-1 \\
c(i_2+1,\ldots, i_s+1,0)  & \mbox{ if } i_1=m-1 \end{array}\right.
\end{equation}
If $s=r$, then
\begin{equation}
  \label{equation:coefficient:determinant}
  c(i_1,i_2,\ldots , i_{r})  = \det(\bv_{i_1}, \bv_{i_2},\ldots,
  \bv_{i_{r}})
\end{equation}
where $\bv_i=(\mu_1^{2^i}, \mu_2^{2^i}, \ldots, \mu_{r}^{2^i})^{\textsf T}$.
We remark that since $\mu_j^{2^m}=\mu_j$, we have $\bv_m=\bv_0$, and we will
read the indices of $\bv_i$ modulo $m$.

Since the $\mu_i$ are linearly independent over $\F_2$,
$c(r-1,r-2,\ldots, 1,0)=\det(\bv_0, \bv_1, \ldots ,\bv_{r-1})$ is
nonzero. For a proof of this fact, see \cite[p.~5]{Gos96} or
\cite[Lemma~3.5]{LidNie97}. Indeed, $\det(\bv_0, \bv_1,
\ldots,\bv_{r-1})$ is usually called a {\it Moore determinant},
which can be viewed as a $q$-analogue of the familiar Vandermonde
determinants. It follows from (\ref{equation:coefficient:shift})
that $\det(\bv_i, \bv_{i+1}, \ldots ,\bv_{i+r-1})=c(i+r-1,\ldots,
i+1,i)\neq 0$ for all $i$. Therefore, any $r$ consecutive vectors
$\bv_{i}, \bv_{i+1}, \ldots ,\bv_{i+r-1}$ from $\bv_0, \bv_1,
\ldots ,\bv_{m-1}$ are linearly independent over $\F_{2^m}$.

The following lemma reveals more surprising relations among the
coefficients of $S(x)$. We will use this lemma in the proof of
Theorem~\ref{theorem:p1}.

\begin{lemma}
  \label{lemma:coefficientproduct}
  $c(r,r-1,\ldots, 2,0)=c(r-1,r-2,\ldots, 1,0)\cdot c(r-1,r-2,\ldots, 2,1)$.
\end{lemma}
\begin{proof} First note that $c(r-1,r-2,\ldots ,1,0)=\det(\bv_0,\bv_1,\ldots ,\bv_{r-1})\neq
0$. In order to prove the lemma, we show that
$$c(r-1,r-2,\ldots ,2,1)=\frac {c(r,r-1,\ldots, 2,0)} {c(r-1,r-2,\ldots,
1,0)}.$$ Now notice that $c(r,r-1,\ldots,
2,0)=\det(\bv_0,\bv_2,\bv_3,\ldots ,\bv_r)$, so we are trying to
prove that $c(r-1,r-2,\ldots ,2,1)$ is a quotient of two
determinants. This motivates us to consider the following linear
system.
\begin{align}
  \label{equation:matrix:1}
  \begin{pmatrix}
    \mu_1 & \mu_1^2 & \cdots & \mu_1^{2^{r-1}} \\
    \mu_2 & \mu_2^2 & & \mu_2^{2^{r-1}} \\
    \vdots & & \ddots & \vdots \\
    \mu_{r} & \mu_{r}^2 & \cdots & \mu_{r}^{2^{r-1}}
  \end{pmatrix}
  \begin{pmatrix}
    b_0 \\
    b_1 \\
    \vdots \\
    b_{r-1}
  \end{pmatrix}
  & =
  \begin{pmatrix}
    \mu_1^{2^{r}} \\
    \mu_2^{2^{r}} \\
    \vdots \\
    \mu_{r}^{2^{r}}
  \end{pmatrix}
\end{align}
The determinant of the coefficient matrix of this system is
$c(r-1,r-2,\ldots ,1,0)\neq 0$. Thus the system has a unique
solution. In particular, by Cramer's rule,
\[ b_1 =
  \frac{\abs{\begin{array}{ccccc}
        \mu_1 & \mu_1^{2^{r}} & \mu_1^{2^2} & \cdots & \mu_1^{2^{r-1}} \\
        \mu_2 & \mu_2^{2^{r}} & \mu_2^{2^2} &  & \mu_2^{2^{r-1}} \\
        \vdots & & & \ddots & \vdots \\
        \mu_{r} & \mu_{r}^{2^{r}} & \mu_{r}^{2^2} & \cdots &
        \mu_{r}^{2^{r-1}}
      \end{array}}}{\abs{\begin{array}{ccccc}
        \mu_1 & \mu_1^2 & \mu_1^{2^2} & \cdots & \mu_1^{2^{r-1}} \\
        \mu_2 & \mu_2^2 & \mu_2^{2^2} &  & \mu_2^{2^{r-1}} \\
        \vdots & \phantom{\mu_{r}^{2^{r}}} & & \ddots & \vdots \\
        \mu_{r} & \mu_{r}^2 & \mu_{r}^{2^2} & \cdots & \mu_{r}^{2^{r-1}}
      \end{array}}}= \frac{c(r,r-1,\ldots , 2, 0)}{c(r-1,\ldots, 1, 0)}. \]

Next we calculate $b_j$'s explicitly in a different way. In
particular, we will show that $b_1=c(r-1,r-2,\ldots ,2,1)$. To
this end, we consider the formal power series
\[
f_t (x) = \left( \sum_{j=0}^{\infty} \mu_t^{2^j} x^{2^j}\right)
\prod_{i=1}^{r} \left(1+\sum_{j=0}^{\infty} \mu_i^{2^j}
x^{2^j}\right)\in\F_{2^m}[[x]]
\]
for $1\le t\le r$. We have
\begin{align*}
  \left( \sum_{j= 0}^{\infty} \mu_t^{2^j} x^{2^j}\right) \prod_{i=1}^{r}
  \left(1+\sum_{j= 0}^{\infty} \mu_i^{2^j} x^{2^j}\right) & = \left( \sum_{j= 0}^{\infty}
    \mu_t^{2^j} x^{2^j}\right)
  \left( 1+\sum_{j= 0}^{\infty} \mu_t^{2^j} x^{2^j}\right) \\
  & \phantom{=\ }\cdot \prod_{\ontop{i=1}{i\neq t}}^{r} \left(1+\sum_{j= 0}^{\infty} \mu_i^{2^j}
    x^{2^j}\right)\\
  & = \mu_t x \cdot \prod_{\ontop{i=1}{i\neq t}}^{r} \left(1+\sum_{j= 0}^{\infty}
    \mu_i^{2^j} x^{2^j}\right)
\end{align*}
For any integer $s\leq r$ and $i_1>i_2>\ldots >i_s\geq 0$, we
denote the coefficient of $x^{2^{i_1}+2^{i_2}+\cdots + 2^{i_s}}$
in $\prod_{i=1}^{r} (1+\sum_{j= 0}^{\infty} \mu_i^{2^j} x^{2^j})$
by $c'(i_1,i_2,\ldots, i_s)$. Note that $c'(i_1,i_2,\ldots, i_s)$
is not necessarily the same as $c(i_1,i_2,\ldots, i_s)$ defined
earlier. The former is the coefficient in a formal power series
$\prod_{i=1}^{r}(1+\sum_{j= 0}^{\infty} \mu_i^{2^j}
x^{2^j})\in\F_{2^m}[[x]]$ while the latter is the coefficient in
$S(x)\in\F_{2^m}[x]/(x^{2^m}-x)$.

Clearly, the coefficient of $x^{2^{r}-1}$ in $\prod_{i=1,i\neq t}^{r}
(1+\sum_{j= 0}^{\infty} \mu_i^{2^j} x^{2^j})$ is 0. This shows that the
coefficient of $x^{2^{r}}$ in $f_t(x)$ is 0. On the other hand, from the
definition of $f_t(x)$, we see that this coefficient is
$\mu_t^{2^{r}}+\sum_{j=0}^{r-1}\mu_t^{2^j} c'(r-1, \ldots, j)$. Thus we
obtain
\begin{align}
  \label{equation:matrix:2}
  \mu_t^{2^{r}} & = \sum_{j=0}^{r-1}\mu_t^{2^j} c'(r-1, \ldots, j+1, j)
\end{align}
for all $1\le t\le r$. Combining (\ref{equation:matrix:1}) and
(\ref{equation:matrix:2}) we have $b_j=c'(r-1, \ldots, j+1, j)$.
In particular, $b_1=c'(r-1,\ldots, 2, 1)$. To finish the proof, we
have to show that $c'(r-1,\ldots, 2, 1)=c(r-1,\ldots, 2, 1)$.
Clearly, it suffices to show that if $s,j_1,\ldots, j_s$ are
integers with $s\le r$ and $0\le j_1, \ldots, j_s\le m-1$ such
that
\begin{equation}\label{congr}
2^{j_1} + 2^{j_2} + \cdots + 2^{j_s} \equiv 2^{r-1} + 2^{r-2} + \cdots + 2
\pmod{2^m - 1}
\end{equation}
then $2^{j_1} + 2^{j_2} + \cdots + 2^{j_s} < 2^m-1$.

For any integer $a$ not divisible by $2^m-1$, we use $w(a)$ to
denote the sum of the digits of $a$ (mod $2^m-1$) written in base
2 representation. Note that if $a+b\not\equiv 0$ (mod $2^m-1$),
then $w(a+b)\leq w(a)+w(b)$, and $w(a)+w(b)-w(a+b)$ is the number
of carries occurred in the addition of $a$ and $b$. Applying this
to the above congruence we see that $s\ge r-1$, thus $s=r$ or
$s=r-1$. Moreover, if $s=r$, then exactly one carry occurs in the
(modular) addition $2^{j_1}+ 2^{j_2} + \cdots + 2^{j_s}$, and if
$s=r-1$, then necessarily $\{j_1,j_2,\ldots ,j_s\}=\{1,2,\ldots
,r-1\}$ and $2^{j_1}+ 2^{j_2} + \cdots + 2^{j_s}<2^m-1$.

Now suppose that $2^{j_1} + 2^{j_2} + \cdots + 2^{j_s} \ge 2^m-1$.
Then, by our previous observation, $s=r$ and exactly one carry
occurs in the addition $2^{j_1}+ 2^{j_2} + \cdots + 2^{j_s}$. This
shows that exactly two or exactly three exponents among
$j_1,j_2,\ldots ,j_s$ are equal. Without loss of generality, we
assume that either $j_1=j_2$ ($j_3>j_4>\cdots >j_s$ and they are
not equal to $j_1$) or $j_1=j_2=j_3$ ($j_4>j_5>\cdots >j_s$ and
they are not equal to $j_1$). In the former case we must have
$m-1=j_1=j_2>j_3>j_4>\cdots >j_s>0$, and
$$2^{j_1} + 2^{j_2} + \cdots
+ 2^{j_s}\equiv 2^{j_3} + 2^{j_4} + \cdots + 2^{j_s}+2^0 \pmod
{2^m-1},$$ contradicting (\ref{congr}). In the latter case, we
must have $m-1=j_1=j_2=j_3>j_4>j_5>\cdots >j_s>0$, and
$$2^{j_1} + 2^{j_2} + \cdots
+ 2^{j_s}\equiv 2^{m-1}+2^{j_4} + 2^{j_5} + \cdots + 2^{j_s}+2^0
\pmod {2^m-1},$$ again contradicting (\ref{congr}). This completes
the proof of the lemma.
\end{proof}

We will also need the following lemma in the proof of
Theorem~\ref{theorem:p1}.

\begin{lemma}
\label{lemma:determinant}
Let
\begin{align*}
  \Delta & = \phantom{+}
  c(m-1,m-2,\ldots, m-r+1, m-r-1) \cdot c(m-2,m-3,\ldots, m-r, 0)\\
  & \phantom{=\,} + c(m-2,m-3,\ldots, m-r+1, m-r-1,0)\cdot
  c(m-1,m-2,\ldots,m-r).
\end{align*}
Then $\Delta\neq 0$.
\end{lemma}
\begin{proof}
Recall that
\begin{align*}
  c(m-1,m-2,\ldots,m-r+1,m-r-1) & = \det (\bv_{m-1}, \bv_{m-2}, \ldots ,
  \bv_{m-r+1} , \bv_{m-r-1})\\
  c(m-2,m-3,\ldots,m-r,0) & = \det (\bv_{m-2}, \bv_{m-3}, \ldots ,
  \bv_{m-r} , \bv_0)\\
  c(m-2,m-3,\ldots,m-r+1,m-r-1,0) & = \det (\bv_{m-2}, \bv_{m-3}, \ldots ,
  \bv_{m-r+1} , \bv_{m-r-1}, \bv_0)\\
  c(m-1,m-2,\ldots,m-r) & =
  \det (\bv_{m-1}, \bv_{m-2}, \ldots , \bv_{m-r})
\end{align*}
Since $\bv_{m-2}, \bv_{m-3}, \ldots , \bv_{m-r-1}$ form a basis of
the $\F_{2^m}$-span of $\{\bv_0, \bv_1,\ldots, \bv_{m-1}\}$, there
exist $\alpha_i$'s and $\beta_i$'s in $\F_{2^m}$ such that
\begin{align}
  \label{equation:determinant:2}
  \bv_{m-1} & = \alpha_{m-2} \bv_{m-2} + \cdots +
  \alpha_{m-r} \bv_{m-r} + \alpha_{m-r-1} \bv_{m-r-1}\\
  \label{equation:determinant:3}
  \bv_0 & = \beta_{m-2} \bv_{m-2} + \cdots +
  \beta_{m-r} \bv_{m-r} + \beta_{m-r-1} \bv_{m-r-1}
\end{align}
Then
\begin{align*}
  c(m-1,m-2,\ldots,m-r+1,m-r-1) & = \alpha_{m-r} \det (\bv_{m-2},\bv_{m-3},
  \ldots ,\bv_{m-r},\bv_{m-r-1})\\
  c(m-2,m-3,\ldots,m-r,0) & = \beta_{m-r-1}\det (\bv_{m-2},\bv_{m-3},
  \ldots ,\bv_{m-r},\bv_{m-r-1})\\
  c(m-2,m-3,\ldots,m-r+1,m-r-1,0) & = \beta_{m-r}\det (\bv_{m-2},\bv_{m-3},
  \ldots ,\bv_{m-r},\bv_{m-r-1})\\
  c(m-1,m-2,\ldots,m-r) & = \alpha_{m-r-1} \det (\bv_{m-2},\bv_{m-3},
  \ldots ,\bv_{m-r},\bv_{m-r-1})
\end{align*}
Hence we have
\[
\Delta=\det (\bv_{m-2},\bv_{m-3}, \ldots,\bv_{m-r},\bv_{m-r-1})^2
\abs{\begin{array}{cc}
\alpha_{m-r} & \alpha_{m-r-1}\\
\beta_{m-r} & \beta_{m-r-1}
\end{array}}
\]
Since $\bv_{m-2}, \bv_{m-3}, \ldots, \bv_{m-r-1}$ are linearly
independent over $\F_{2^m}$, $\det (\bv_{m-2}, \bv_{m-3}, \ldots,
\bv_{m-r-1})$ is nonzero. The second determinant in the right hand
side (RHS) of the above equation has to be nonzero for otherwise
(\ref{equation:determinant:2}) and (\ref{equation:determinant:3})
give a dependence relation for the $r$ consecutive vectors
$\bv_{m-1}, \bv_{m-2},\ldots, \bv_{m-r+1}$, $\bv_0$ (note that
$\bv_0=\bv_m$). This shows that $\Delta\neq 0$.
\end{proof}

We are now ready to give the proof of Theorem~\ref{theorem:p1}.

\begin{proof}[Proof of Theorem~\ref{theorem:p1}] Recall that we assume
the defining equation for $A$ is
  \[
  \prod_{i=1}^r (1+\trace(\mu_i x))=1,
  \]
  where $r=m-d$. Suppose that $\trace(a_0)=0$. Then
  $(1+\trace(\sum_{i=1}^{d-1}a_i\lambda^{2^i-1}))=0$ for all $\lambda\in
  A\setminus \{0\}$. Thus, the function from $\F_{2^m}$ to $\F_{2^m}$
  associated with the polynomial
  $(1+\trace(\sum_{i=1}^{d-1}a_i x^{2^i-1})) \prod_{i=1}^{r}
  (1+\trace(\mu_i x))$ is the characteristic function of $\{0\}$ in $\F_{2^m}$.  Hence, we have
  \begin{align}
    \label{equation:tr:0}
    \left(1+\trace\left( \sum_{j=1}^{d-1} a_j x^{2^j-1}\right)\right)
    \prod_{i=1}^{r} (1+\trace(\mu_i x)) & \equiv x^{2^m-1}-1\pmod{x^{2^m}-x}.
  \end{align}
  The binary representation of the exponent of $x^{2^m-1}$ (in the LHS of
  (\ref{equation:tr:0})) is $11\ldots 1$ ($m$ ones altogether).
  Throughout this paper we write the most significant bit (i.e., the
  $(m-1)$th bit) to the least significant bit (i.e., the $0$th bit) from
  left to right. Note that the binary representation of the exponent of any
  term in $\prod_{i=1}^{r} (1+\trace(\mu_i x))$ cannot have more than $r$
  ones. The binary representation of the exponent of any term in
  $(1+\trace(\sum_{j=1}^{d-1} a_j x^{2^j-1}))$ has at most $d-1$ ones. Thus,
  the maximum number of ones in the binary representation of the exponent of
  any term on the left hand side of (\ref{equation:tr:0}) is $r+(d-1)=m-1$.
  Therefore the coefficient of $x^{2^m-1}$ on the LHS of
  (\ref{equation:tr:0})  is 0. This contradicts (\ref{equation:tr:0}). So
  this case does not occur.

  From now on we assume that $\trace(a_0)=1$. Then
  $\trace(\sum_{i=1}^{d-1}a_i\lambda^{2^i-1})=0$ for all $\lambda\in
  A\setminus \{0\}$. Therefore the function from $\F_{2^m}$ to $\F_{2^m}$
  associated with the polynomial $$\trace(\sum_{j=1}^{d-1} a_j x^{2^j-1})
  \prod_{i=1}^{r} (1+\trace(\mu_i x))\in \F_{2^m}[x]$$ is the zero function.
  That is, in $\F_{2^m}[x]$, we have the following congruence.
\begin{align}
  \label{equation:traceproduct}
  \trace\left( \sum_{j=1}^{d-1} a_j x^{2^j-1}\right) \prod_{i=1}^{r}
  (1+\trace(\mu_i x)) & \equiv 0\pmod{x^{2^m}-x}
\end{align}
For later use, we let $T(x)$ and $S(x)$ be polynomials in
$\F_{2^m}[x]$ of degree less than or equal to $2^m-1$ such that
$T(x)\equiv \trace( \sum_{j=1}^{d-1} a_j x^{2^j-1})\pmod
{x^{2^m}-x}$ and $S(x)\equiv \prod_{i=1}^{r} (1+\trace(\mu_i
x))\pmod {x^{2^m}-x}$.

Now the proof proceeds as follows. We will first prove that
$a_{d-1}=a_{d-2}=0$. Next we will show that the ``upper half'' 
coefficients of $p(x)$ are zero. More precisely, we prove that
$a_{m-d+1}=a_{m-d+2}=\cdots =a_{d-3}=0$. Finally we show that the
``lower half'' coefficients of $p(x)$ are also zero. That is,
$a_2=a_3=\cdots =a_{m-d}=0$ (here we assume that $m-d\geq 2$).

{\bfseries Claim: \boldmath$a_{d-1}=a_{d-2}=0$}. Consider the
coefficient of the monomial $x^{(2^{m-1}-1)-2^{d-2}}$ in
$T(x)\cdot S(x)$, i.e., the left hand side (LHS) of
(\ref{equation:traceproduct}). The binary expansion of its
exponent is
\[
  0\overbrace{1\ldots 11}^{r}0\overbrace{1\ldots 1}^{d-2}.
\]
  The number of 1's in this expansion is $(m-2)$. The maximum number
  of 1's in the exponent of
  any summand in $S(x)$ is $r$ and the maximum number of 1's in the
  exponent of any summand in $T(x)$ is $d-1$. When adding two exponents
  (written in their binary representations), any carry that may occur reduces
  the number of 1's in the sum. Since we are interested in an
  exponent whose number of 1's is $(m-2)$, it can only be obtained as a sum
  of two exponents (one is the exponent of a summand in $T(x)$, the other
  in $S(x)$) with at most one carry.

If $(2^{m-1}-1)-2^{d-2}$ is obtained as a sum without carry then
there is only one possibility.
\begin{align*}
  0\overbrace{1\ldots 11}^{r}0\overbrace{1\ldots 1}^{d-2} & =
  0\overbrace{1\ldots 11}^{r}000\ldots 00 + 00\ldots 000\overbrace{11\ldots 11}^{d-2}
\end{align*}
Using the assumption that $2d>m+2$, we see that $r<d-2$ and thus,
the $d-2$ consecutive 1's have to come from the term
$x^{2^{d-2}-1}$ in $T(x)$, whose coefficient is $a_{d-2}$.

If $(2^{m-1}-1)-2^{d-2}$ is obtained as a sum with exactly one
carry, then that carry has to happen at position $d-2$ and so
\begin{align*}
0\overbrace{1\ldots 11}^{r}0\overbrace{1\ldots
  1}^{d-2} & = 0\overbrace{1\ldots 1}^{r}0100\ldots 00 + 00\ldots 00\overbrace{111\ldots 11}^{d-1}
\end{align*}
Again, the $d-1$ consecutive 1's have to come from the term
$x^{2^{d-1}-1}$ in $T(x)$, whose coefficient is $a_{d-1}$. Hence
by (\ref{equation:traceproduct}), we have
\begin{align}
  & \phantom{{}+{}} c(m-2, m-3, \ldots, d, d-1)\cdot a_{d-2}\nonumber\\
  & {}+ c(m-2,m-3,\ldots, d, d-2)\cdot a_{d-1}\label{equation:coefficient:1}\\
  & = 0.\nonumber
\end{align}

Next we look at the coefficient of $x^{(2^{m-1}-1)-2^{d-1}}$ in
$T(x)\cdot S(x)$. As before, the number of 1's in the binary
expansion of $(2^{m-1}-1)-2^{d-1}$ is $m-2$. Hence at most one
carry may occur. Again, using $r-1<d-1$ there are only three ways
of obtaining $(2^{m-1}-1)-2^{d-1}$ as a sum of two exponents
without carry.
\begin{align*}
  0\overbrace{1\ldots 1}^{r-1}0\overbrace{11\ldots 1}^{d-1} & =
  0\overbrace{1\ldots 11}^{r-1}000\ldots 00 +
  00\ldots 000\overbrace{11\ldots 11}^{d-1}\\
  & = 0\overbrace{1\ldots 11}^{r-1}01\overbrace{0\ldots 00}^{d-2} + 00\ldots 0000\overbrace{1\ldots 11}^{d-2}\\
  & = 0\overbrace{1\ldots 11}^{r-1}0\overbrace{00\ldots 0}^{d-2}1 + 00\ldots 000\overbrace{11\ldots 1}^{d-2}0 \intertext{If a
    carry occurs, then it has to be at position $d-1$.}  0\overbrace{1\ldots 1}^{r-1}0\overbrace{11\ldots 1}^{d-1} &
  = 0\overbrace{1\ldots 1}^{r-2}0100\ldots 01 + 00\ldots 00\overbrace{111\ldots 1}^{d-1}0
\end{align*}
It follows from (\ref{equation:traceproduct}) that
\begin{align}
  & \phantom{{}+{}} c(m-2, m-3, \ldots, d)\cdot a_{d-1}\nonumber\\
  & {}+ c(m-2,m-3,\ldots, d, d-2)\cdot a_{d-2}\nonumber\\
  \label{equation:coefficient:2}
  & {}+ c(m-2, m-3, \ldots, d,0) \cdot a_{d-2}^2\\
  & {}+ c(m-2, m-3, \ldots, d+1,d-1,0) \cdot a_{d-1}^2\nonumber\\
  & = 0.\nonumber
  \intertext{Now we claim that}
  & \phantom{{}+{}} c(m-2, m-3, \ldots, d)\cdot a_{d-1}\nonumber\\
  & {}+ c(m-2,m-3,\ldots, d, d-2)\cdot a_{d-2}\label{equation:coefficient:3}\\
  & = 0\nonumber.
\end{align}
In order to prove (\ref{equation:coefficient:3}), we will show
that
\begin{align*}
  \abs{\begin{array}{cc}
      c(m-2, m-3, \ldots, d,d-2) & c(m-2, m-3,\ldots, d-1)\\
      c(m-2,m-3,\ldots, d) & c(m-2,m-3,\ldots, d, d-2)
  \end{array}} & = 0
\end{align*}
Once we prove this, it is clear that
(\ref{equation:coefficient:3}) will follow from
(\ref{equation:coefficient:1}). Hence we need to show that
\begin{align}
  c(m-2, m-3, \ldots, d,d-2)^2 & = c(m-2, m-3,\ldots, d-1)\nonumber\\
  \label{equation:determinant:1}
  & \phantom{=\ } \cdot c(m-2,m-3,\ldots, d) \intertext{which, by
  (\ref{equation:coefficient:shift}) is the same as}
  c(m-1, m-2, \ldots, d+1,d-1) & = c(m-2, m-3,\ldots, d-1)\nonumber\\
  & \phantom{=\ } \cdot c(m-2,m-3,\ldots, d)\nonumber
  \intertext{Making appropriate shifts using (\ref{equation:coefficient:shift}),
  the above equation is further equivalent to}
  c(r,r-1 \ldots, 2,0)& ={c(r-1,r-2 \ldots, 1, 0)}\cdot c(r-1,\ldots, 2, 1).\nonumber
\end{align}
Hence, by Lemma~\ref{lemma:coefficientproduct}, we have proved
(\ref{equation:coefficient:3}).

Now the combination of (\ref{equation:coefficient:1}),
(\ref{equation:coefficient:2}), and (\ref{equation:coefficient:3})
yields that
\begin{align}
\label{equation:matrix:3}
\begin{pmatrix}
c(m-2, \ldots, d, d-2)^2 & c(m-2, \ldots, d, d-1)^2\\
c(m-2, \ldots, d+1, d-1, 0) & c(m-2, \ldots, d,0)
\end{pmatrix}
\begin{pmatrix}
a_{d-1}^2\\
a_{d-2}^2
\end{pmatrix}
& =
\begin{pmatrix}
0\\
0
\end{pmatrix}
\end{align}
By Lemma~\ref{lemma:determinant} the determinant of the coefficient matrix
in (\ref{equation:matrix:3}) is nonzero and thus, $a_{d-1}=a_{d-2}=0$.

{\bfseries Claim: \boldmath$a_{d-3}=\cdots =a_{r+1}=0$}. Now let $d-2>k>r$
and suppose that $a_j=0$ for all $d-1\ge j>k$. We want to show that $a_k=0$.
To this end, consider the coefficient of $x^{(2^{m-1}-2^{d-1})+(2^k-1)}$ in
$T(x)\cdot S(x)$. Since $k>r$ there is only one way of attaining this
exponent when multiplying $T(x)$ and $S(x)$.
\begin{align*}
0\overbrace{1\ldots 1}^{r}0\ldots 0\overbrace{11\ldots 1}^{k} & =
0\overbrace{1\ldots 11}^{r}000\ldots 00 + 00\ldots 000\overbrace{11\ldots 11}^{k}
\end{align*}
Hence by (\ref{equation:traceproduct}), $c(m-2, m-3, \ldots,
d-1)\cdot a_k=0$. As noted before, $c(m-2, m-3, \ldots, d-1)\neq
0$ so we have $a_k=0$.

At this point we note that if $d=m-1$, i.e., $r=1$, then the above
two claims already show that $a_2=a_3=\cdots =a_{d-1}=0$, and the
theorem is proved in this case. So from now on, we assume that
$m-1>d>\frac{m}{2}+1$. Also we will assume that $m\geq 10$. The case
where $5\leq m\leq 8$ will be dealt with separately at the very
end of the proof.

{\bfseries Claim: \boldmath$a_3=\cdots =a_{r}=0$}. For any integer
$t$, $3\leq t\leq r$, suppose that $a_j=0$ for all $j>t$, we will
prove that $a_t=0$. Here we need the following result, whose proof
will be given right after our proof of Theorem~\ref{theorem:p1}.

{\bfseries Result 1:} Assume that $m\geq 10$ and
  $\lfloor\frac{m-3}{2}\rfloor\ge r\ge t\ge 3$. There exist $0=i_1<\cdots<i_{r}\le m-t-3$ such that
\begin{enumerate}
  \renewcommand{\labelenumi}{\emph{(\roman{enumi})}}
\item $c(i_1,i_2,\ldots ,i_{r})\neq 0$, and
\item the number of consecutive integers in the set $\{i_1,i_2,\ldots
  ,i_{r}\}$ is less than or equal to $t-1$.
\end{enumerate}

With Result 1, we will look at the coefficient of $x^{(2^{m-1}-2^{m-t-1}) +
\sum_{j=1}^{r}2^{i_j}}$ in $T(x)\cdot S(x)$, i.e., the LHS of
(\ref{equation:traceproduct}). Note that the exponent of this
monomial has the $m$-bit binary representation
\[
  0\underbrace{11\ldots 1}_{t}0\underbrace{0\ldots 1\ldots 1\ldots 1}_{m-t-2},
\]
where at $i_j$-th bit, there is a 1, for each $j=1,2,\ldots ,r$.

Since the number of consecutive integers in the set
$\{i_1,i_2,\ldots ,i_r\}$ is less than or equal to $t-1$, there is
only one way to get the term $x^{(2^{m-1}-2^{m-t-1}) +
\sum_{j=1}^{r}2^{i_j}}$ when multiplying $T(x)$ with $S(x)$,
namely
\[
0\underbrace{11\ldots 1}_{t}0\underbrace{0\ldots 1\ldots 1\ldots
  1}_{m-t-2}=0\underbrace{00\ldots 0}_{t}0\underbrace{0\ldots 1\ldots
  1\ldots 1}_{m-t-2}+0\underbrace{11\ldots 1}_{t}0\underbrace{00\ldots
  0}_{m-t-2}.
\]
Therefore, the coefficient of $x^{(2^{m-1}-2^{m-t-1}) +
\sum_{j=1}^{r}
  2^{i_j}}$ in $T(x)\cdot S(x)$ is
$c(i_1,i_2,\ldots ,i_{r})\cdot a_t^{2^{m-t-1}}$. It follows now
from (\ref{equation:traceproduct}) that
\begin{align*}
  c(i_1,i_2,\ldots ,i_{r})\cdot a_t^{2^{m-t-1}} & = 0.
\end{align*}
Noting that $c(i_1,i_2,\ldots ,i_{r})\neq 0$ we have $a_t=0$.

{\bfseries Claim: \boldmath$a_2=0$}. Suppose that $a_2\neq 0$. Let
$Q (x)=\trace(a_2 x^3+a_1 x)$ and let $V=\F_{2^m}$. Note that
since $\trace(a_0)=1$, the assumption that $\trace(p(\lambda))=1$
for all $\lambda\in A\setminus\{ 0\}$ implies that $Q(\lambda)=0$
for all $\lambda\in A$, where $\abs{A}=2^d$. The map $Q:V\to\F_2$
is a quadratic form with associated bilinear form
\begin{align*}
  B(x,y) & = Q (x+y) - Q (x) - Q (y)\\
  & = \trace(a_2(xy^2+yx^2)).
\end{align*}
We will show that the maximum dimension of a subspace of $V$ on
which $Q$ vanishes is less than $d$. This will force $a_2=0$.

Let $\rad V=\{x\in V\mid B(x,y)=0,\,\forall {y\in V}\}$. Note that
in even characteristic $Q$ does not have to be zero on $\rad V$.
Therefore we consider $V_0=\{x\in\rad V\mid Q (x)=0\}$. We call
$Q$ {\it nonsingular} if $V_0=\{0\}$. By Witt's theorem, the
maximum dimension of a totally singular subspace of a nonsingular
quadratic space $(V,Q)$ is at most $\floor{\frac{1}{2}\dim V}$. In
our case we have $\rad V=\{x\in V\mid x=a_2 x^4\}$. In particular,
$\dim V_0\le 2$. If $Q$ is nonsingular then the maximum dimension
of a totally singular subspace is at most $\floor{\frac{m}{2}}$.
If $Q$ is singular then we consider the induced (nonsingular)
quadratic form $\bar Q:V/V_0\to\F_2$. The maximum dimension of a
subspace $U$ of $V/V_0$ on which $\bar Q$ vanishes is at most
$\frac{1}{2}(m-\dim V_0)$. The maximum dimension of a subspace of
$V$ on which $Q$ vanishes is less than or equal to $\dim(U\perp
V_0)\leq \frac{1}{2}(m+\dim V_0)\leq \frac {m}{2}+1$. It follows
that in either case the maximum dimension of a subspace of $V$
on which $Q$ vanishes is less than $d$, hence $a_2$ has to be 0.

Finally we deal with the case where $5\leq m\leq 8$. When $m=5$ or
6, there is no admissible $d$ satisfying the restriction that $m-1>d>\frac
{m}{2}+1$. When $m=7$ (resp. $m=8$), the only admissible $d$ is 5
(resp. 6). In both cases, $r=m-d=2$, and by the first two claims,
we have $a_3=a_4=\cdots =a_{d-1}=0$. Now by the same argument
using quadratic forms as above, we can further prove that $a_2=0$.

The proof of the theorem will be complete once we proof Result 1
above.
\end{proof}

We now give the promised proof of Result 1. This result can be
thought as a generalization of the fact that a Moore determinant
is nonzero, and it may be of independent interest. The proof of
Result 1 we give here is elementary, but quite technical. The
reader may want to skip the proof in a first reading of the paper.

We state Result 1 formally as
\begin{theorem}
  \label{lemma:gap}
  Let $m, r, t$ be positive integers, and let $\mu_1, \ldots ,
  \mu_r\in\F_{2^m}$ be linearly independent over $\F_2$. If
  $m\geq 10$ and $\lfloor\frac{m-3}{2}\rfloor\ge r\ge t\ge 3$,
  then there exist $0=i_1<i_2<\cdots <i_r\leq m-(t+3)$ such that
  \begin{enumerate}
  \item $\det
    \begin{pmatrix}
      \mu_1^{2^{i_1}} & \mu_2^{2^{i_1}} & \cdots & \mu_r^{2^{i_1}}\\
      \vdots & & \ddots & \vdots\\
      \mu_1^{2^{i_r}} & \mu_2^{2^{i_r}} & \cdots & \mu_r^{2^{i_r}}\\
    \end{pmatrix}\neq 0$, and
  \item the number of consecutive integers in the set
    $\{i_1,i_2,\ldots ,i_r\}$ is at most $t-1$.
  \end{enumerate}
\end{theorem}
We first fix some notation. Let $V$ be the $\F_{2^m}$-span of
$\{\bv_0, \ldots , \bv_{m-1}\}$, where
$\bv_i=(\mu_1^{2^i},\mu_2^{2^i}, \ldots, \mu_r^{2^i})^{\textsf
T}$. As before, all indices of the vectors $\bv_i$ are to be read
modulo $m$. We have $\dim_{\F_{2^m}} V=r$ and $\{\bv_i,\bv_{i+1},
\ldots ,\bv_{i+r-1}\}$ is a basis of $V$ for all $i\geq 0$
\cite[Lemma~3.5]{LidNie97}. By $\bv_i^{2^j}$ we mean
component-wise exponentiation of $\bv_i$ by $2^j$. Hence
$\bv_i^{2^j}=\bv_{i+j}$. We will use binary vectors to denote
subsets of $\{\bv_0, \ldots , \bv_{m-1}\}$ as follows.  Let
$\bu=(u_0,u_1,\ldots, u_i)\in\F_2^{i+1}$ be a vector of length
$i+1$.  By $\Lambda(\bu)$ we will denote the $\F_{2^m}$-span of
$\{ \bv_j \mid u_j=1\}$.  We also allow concatenation of binary
vectors. If $\bu=(u_0,u_1,\ldots, u_i)$ and
$\bu'=(u_0',u_1',\ldots, u_j')$ then
\[
\bu * \bu' = (u_0,u_1,\ldots, u_i,u_0', u_1',\ldots, u_j').
\]
If we concatenate several copies, say $i\ge 1$, of the same vector $\bu$
then we denote the resulting vector by $\bu^{*i}$. Sometimes it may happen
that we have a concatenated vector $\bu'*\bu^{*i}$ with $i=0$. In this case
we assume that no copy of $\bu$ had been appended to $\bu'$, that is,
$\bu'*\bu^{*0}=\bu'$.

Now Theorem~\ref{lemma:gap} can be reformulated as follows.
\begin{theorem*}[\protect{\bfseries\ref{lemma:gap}'}]
  For $r\le\floor{\frac{m-3}{2}}$ there exists a binary vector\/ $\bw$ of
  length at most $m-(t+2)$ such that $\Lambda(\bw)=V$ and the number of
  consecutive 1's in\/ $\bw$ is at most $t-1$.
\end{theorem*}

  First of all, note that it suffices to prove the theorem in the case
  where $r$ is equal to $\floor{\frac{m-3}{2}}$. Indeed, if we have found a vector $\bw$
  for $\floor{\frac{m-3}{2}}=R$ then the same vector $\bw$ will satisfy our
  requirements for smaller $r$. The reason is as follows. Suppose that $r<R$. We can
  extend the set $\{ \mu_1, \mu_2, \ldots , \mu_r\}$ to a set of $R$
  elements
  $\{ \mu_1, \mu_2, \ldots , \mu_r, \ldots, \mu_R\}$ in $\F_{2^m}$ that are linearly
  independent over $\F_2$. By assumption, we can find $0=i_1<i_2<\cdots
  <i_R<m-(t+3)$ such that $\bv_{i_1}, \bv_{i_2}, \ldots, \bv_{i_R}$ form a
  basis of $\F_{2^m}^R$, where $\bv_{i_j}=(\mu_1^{2^{i_j}},\mu_2^{2^{i_j}},
  \ldots,\mu_R^{2^{i_j}})^{\textsf T}$. Let $\bv_{i_j}'$ be the projection of $\bv_{i_j}$
  onto the first $r$ coordinates, that is,
  $\bv_{i_j}'=(\mu_1^{2^{i_j}}, \mu_2^{2^{i_j}}, \ldots,
  \mu_r^{2^{i_j}})^{\textsf T}$ for $1\le j\le R$.  Then $\{\bv_{i_1}',
  \bv_{i_2}', \ldots, \bv_{i_R}'\}$ spans $\F_{2^m}^r$. Hence this set
  contains $r$ vectors $\bv_{i_{j_1}}', \bv_{i_{j_2}}', \ldots,
  \bv_{i_{j_r}}'$ that are linearly independent over $\F_{2^m}$. By
  assumption, $0\le i_{j_1}<i_{j_2}<\cdots <i_{j_r}<m-(t+3)$ and the number
  of consecutive integers in $\{ i_{j_1}, i_{j_2}, \ldots, i_{j_r}\}$ is at
  most $t-1$. If $i_{j_1}\neq 0$ then it is clear that we can use $\{ 0,
  i_{j_2}-i_{j_1}, \ldots, i_{j_r}-i_{j_1}\}$ instead. From now on, we will
  assume that $r=\floor{\frac{m-3}{2}}$.

  We write $r=kt+a$, where $0\le a\le t-1$. Since $r\geq t$, we have $k\geq 1$.
  Let $\ba=(1,\ldots, 1)\in\F_2^a$, $\bu=(0,1,\ldots,1)\in\F_2^t$, $\bbu=(1,0,\ldots,0)\in\F_2^t$, and
  $\bzero=(0,\ldots,0)\in\F_2^t$. Then $\dim\Lambda(\ba*\bu^{*k})=r-k$ since
  $\ba*\bu^{*k}$ is a vector of length $r$ with exactly $k$ zeros. We will
  append copies of $\bu$ or $\bbu$ to $\ba*\bu^{*k}$ to describe a set of
  vectors $\bv_i$ that generate $V$. Note that by appending $\bu$ to
  $\ba*\bu^{*(k+b)}$, $0\le b<k$, we have
  \[
  \dim\Lambda(\ba * \bu^{*(k+b+1)})\geq \dim\Lambda(\ba * \bu^{*(k+b)}).
  \]
  \begin{lemma}\
    \label{lemma:nestedsequence}
    \begin{enumerate}
    \item \label{lemma:nestedsequence:i} If $\Lambda (\ba*\bu^{*(k+b+1)}) =
      \Lambda (\ba*\bu^{*(k+b)})$, then $\Lambda (\ba*\bu^{*(k+b+i)}) =
      \Lambda (\ba*\bu^{*(k+b)})$ for any positive integer $i$.
    \item \label{lemma:nestedsequence:ii} If $\Lambda (\ba*\bu^{*(k+1)})
      =\Lambda (\ba*\bu^{*k}) + \F_{2^m}\bv_{r+\ell}$ where $1\le \ell\le
      t-1$, then $\Lambda (\ba*\bu^{*(k+b+1)}) =\Lambda (\ba*\bu^{*(k+b)}) +
      \F_{2^m}\bv_{r+bt+\ell}$ for any positive integer $b$.
    \item \label{lemma:nestedsequence:iii} Let $\by$ be a binary vector of length
      $\ell$ and $\bone_t=(1,1,\ldots,1)\in \F_{2}^{t}$. Suppose
      $\Lambda(\by)$ is a proper subspace in $V$ and there exists a vector
      $\bz\in \F_{2}^{t}$ such that $\{\bv_i\mid (\bz*\by)_i=1\}\subseteq
      \{\bv_i\mid (\by*\bone_t)_i=1\}$. Then $\Lambda(\by)\subsetneq
      \Lambda(\by*\bone_t)$.
    \end{enumerate}
  \end{lemma}
  \begin{proof} (1). Observe that since $\Lambda(\ba * \bu^{*(k+b)})=\Lambda(\ba *
      \bu^{*(k+b+1)})$, we have dependence relations
      $\bv_{r+bt+j}=\sum_{i=0}^{r+bt-1} c_i \bv_i$ for $1\le j\le t-1$ where
      $c_i=0$ if $i$ is of the form $a+st$.  This gives dependence relations
      $\bv_{r+bt+j}^{2^t}=\bv_{r+(b+1)t+j}=\sum_{i=0}^{r+bt-1} c_i^{2^t}
      \bv_{i+t}$, hence $\Lambda(\ba * \bu^{*(k+b+2)})\subseteq \Lambda(\ba
      * \bu^{*(k+b+1)})$.

      (2). Observe that our assumption implies that for $1\leq j\leq t-1$
      with $j\neq \ell$, we have dependence relations
      $\bv_{r+j}=c_{r+\ell}\bv_{r+\ell}+ \sum_{i=0}^{r-1} c_i \bv_i$ where
      $c_i=0$ if $i$ is of the form $a+st$. As before, we then obtain the
      relation $\bv_{r+bt+j}=c_{r+\ell}^{2^{bt}}\bv_{r+bt+\ell}+
      \sum_{i=0}^{r-1} c_i^{2^{bt}} \bv_{i+bt}$ where $c_i=0$ if $i$ is of
      the form $a+st$. Clearly, $\sum_{i=0}^{r-1} c_i^{2^{bt}} \bv_{i+bt}\in
      \Lambda(\ba * \bu^{*(k+b)})$ as $c_i=0$ if $i$ is of the form $a+st$.
      We thus obtain (\ref{lemma:nestedsequence:ii}).

     (3). Suppose $\Lambda(\by)= \Lambda(\by* \bone_t)$. Then the $t$
      consecutive vectors $\bv_{\ell},\bv_{\ell +1}, \ldots,
      \bv_{\ell+(t-1)}$ are all in $\Lambda(\by)$.  On the other hand,
      \[ \{x^{2^t}\mid x\in \Lambda(\by)\}\subseteq \Lambda(\bz*\by)\subseteq
      \Lambda(\by* \bone_t)=\Lambda(\by).\]
      It follows that for any $x\in\Lambda(\by)$, we have $x^{2^t}\in
      \Lambda(\by)$.  In particular, $\bv_{\ell+it}$,$\bv_{\ell +1+it}$,
      $\ldots$, $\bv_{\ell+(t-1)+it}$ are all in $\Lambda(\by)$ for every positive
      integer $i$. We thus have $\bv_{0},\ldots, \bv_{r-1}\in \Lambda(\by)$.
      This contradicts our assumption that $V\neq \Lambda(\by)$.
\end{proof}

  We are now ready to prove Theorem~\ref{lemma:gap}'.  Recall that
  $k\geq 1$ and we may assume $r=\floor{\frac{m-3}{2}}$.

\begin{proof} We will consider two cases.

{\bfseries Case \boldmath$\Lambda(\ba * \bu^{*k} * \bu^{*i})\neq
V$ for
  all $i>0$:} In this case the dimensions of the subspaces in
  the nested sequence $$\Lambda (\ba*\bu^{*k})\subseteq
  \Lambda(\ba*\bu^{*(k+1)})\subseteq \cdots
  \Lambda(\ba*\bu^{*(k+i)})\subseteq\cdots$$ stop growing eventually. Let $b$ be the
  largest integer such that $\dim\Lambda(\ba*\bu^{*(k+b)})>\dim\Lambda(\ba
  *\bu^{*(k+b-1)})$. Since $r>\dim\Lambda(\ba*\bu^{*(k+b)})\ge r-k+b$, we
  have $0\le b<k$. By repeated application of
  Lemma~\ref{lemma:nestedsequence}, part~\ref{lemma:nestedsequence:iii}, we see
  that
  \[V =
  \begin{cases}
    \Lambda(\ba * \bu^{*(k+b+1)} * \bone_t^{*(k-b-1)}) & \text{if\ }\dim
    \Lambda(\ba *\bu^{*(k+b)})> r-k+b, \\
    \Lambda(\ba * \bu^{*(k+b)} * \bone_t^{*(k-b)}) & \text{if\ }
    \dim\Lambda(\ba *\bu^{*(k+b)})= r-k+b.
  \end{cases}
  \]
  Since $\Lambda(\ba *\bu^{*(k+b+i)})= \Lambda(\ba*\bu^{*(k+b)})$ for all
  positive integer $i$,
  \[V =
  \begin{cases}
    \Lambda(\ba * \bu^{*(k+b)} *\bzero *\bbu^{*(k-b-1)}) & \text{if\ }
    \dim\Lambda(\ba *\bu^{*(k+b)})> r-k+b, \\
    \Lambda(\ba * \bu^{*(k+b)} * \bbu^{*(k-b)}) & \text{if\ }
    \dim\Lambda(\ba *\bu^{*(k+b)})= r-k+b.
  \end{cases}
  \]

  {\slshape Subcase \boldmath$\dim \Lambda(\ba *\bu^{*(k+b)})> r-k+b$:} We
  define $\bw$ to be the vector obtained by dropping the last $(t-1)$ zeros
  from the last copy of $\bbu$ in the vector $\ba * \bu^{*(k+b)} *\bzero
  *\bbu^{*(k-b-1)}$. Note that the length of $\bw$ is at most $m-(t+2)$,
  $\Lambda(\bw)=V$, and the number of consecutive $1$'s in $\bw$ is at most
  $t-1$.

  {\slshape Subcase \boldmath$\dim \Lambda(\ba *\bu^{*(k+b)})=r-k+b\text{
      and }b>0$:} Appending the $(k+b)$-th copy of $\bu$ increased the
  dimension of $\Lambda(\ba *\bu^{*(k+b-1)})$ by exactly one,
  i.e., $\dim \Lambda(\ba *\bu^{*(k+b)})=1+\dim \Lambda(\ba
  *\bu^{*(k+b-1)})$. Thus, $\Lambda(\ba
  *\bu^{*(k+b)})=\Lambda(\ba *\bu^{*(k+b-1)})+\F_{2^m}\bv_{r+(b-1)t+i}$ for
  some $1\leq i \leq t-1$.  Let $\bu_i\in\F_2^t$ be the vector with
  $(i+1)$-th entry being one and all other entries $0$. Then it is clear that
  $\Lambda(\ba *\bu^{*(k+b)})=\Lambda(\ba *\bu^{*(k+b-1)}*\bu_i)$. Recall
  that $V=\Lambda(\ba * \bu^{*(k+b)} * \bbu^{*(k-b)})$. Therefore, we deduce
  \[V= \Lambda(\ba *\bu^{*(k+b-1)}*\bu_i*\bbu^{*(k-b)}).\]

  To find the required vector $\bw$, we simply drop the last $(t-1)$ zeros from
  the last copy of $\bbu$ in $\ba*\bu^{*(k+b-1)}*\bu_i*\bbu^{*(k-b)}$.
  Clearly, the resulting vector is of length $r+(k-1)t+1$ which is at most
  $m-(t+2)$ and satisfies what we require.

  {\slshape Subcase \boldmath$\dim \Lambda(\ba *\bu^{*(k+b)})=r-k+b\text{
      and }b=0$:} In this case we have
  $\Lambda(\ba*\bu^{*k})=\Lambda(\ba*\bu^{*(k+i)})$ for all $i>0$. It
  follows that for any $0\leq j \leq r-1$ and for any $i$, $\bv_{j+it}\in
  \Lambda(\ba*\bu^{*k})$ if and only if $j\not\equiv a \pmod t$.  As $\bv_{a+m}\notin
  \Lambda(\ba*\bu^{*k})$, it follows that $a+m \equiv a \pmod t$. Hence, we
  have $t\doesdivide m$. Since we may assume $r=\floor{\frac{m-3}{2}}$ it
  follows that $t=2a+3$, $t=2a+4$, or $t=a+2$. In each case, $0\le a\le t-2$.

  First we will assume $t-2\ge a\ge 1$. It follows that $\bv_{t-1}\in
  \Lambda(\ba*\bu^{*k})$ and thus $\bv_{m-1}\in \Lambda(\ba *\bu^{*k})$.
  Recall that any $r$ consecutive vectors in $\{\bv_0,\bv_1,\ldots ,\bv_{m-1}\}$
  are linearly independent. In
  particular, $\{\bv_{m-1}, \bv_0, \ldots, \bv_{r-2}\}$ are linearly
  independent.  Let $\bz=(0,\ldots, 0,1)\in \F_2^{m-r-kt}$. It is clear that
  $\Lambda(\ba*\bu^{k-1}* (0,\underbrace{1,1,\ldots,1}_{t-2},0)
  *\bzero^{*k}*\bz)$ is an $(r-k)$ dimensional subspace in
  $\Lambda(\ba*\bu^{*k})$. As $\dim \Lambda(\ba*\bu^{*k})=r-k$, it
  follows that
  \[
  \Lambda(\ba*\bu^{*(k-1)} * (0,\underbrace{1,1,\ldots,1}_{t-2},0) *
  \bzero^{*k}*\bz)=\Lambda(\ba*\bu^{*k}).
  \]
  Consequently, by
  Lemma~\ref{lemma:nestedsequence}, part~\ref{lemma:nestedsequence:iii}, we
  conclude that
  \[
  \Lambda(\ba*\bu^{*(k-1)} * (0,\underbrace{1,1,\ldots,1}_{t-2},0) *
  \bbu^{*k}*\bz)=\Lambda(\ba*\bu^{*k}*\bbu^{*k}*\bz)=V.
  \]
  The vector $\ba *\bu^{*(k-1)} *(0,\overbrace{1,1,\ldots,1}^{t-2},0)
  *\bbu^{*k} *\bz$ does not have more than $t-1$ consecutive 1's since $a\le
  t-2$. Shifting this vector by one to the right it follows that
  \[
  V=\Lambda((1)*\ba*\bu^{*k-1}*(0,\underbrace{1,1,\ldots,1}_{t-2},0)*\bbu^{*k-1} *(1)).
  \]
  The length of the vector
  $\bw=(1)*\ba*\bu^{*k-1}*(0,\overbrace{1,1,\ldots,1}^{t-2},0)*\bbu^{*k-1}
  *(1)$ is $r+(k-1)t+2$, which is at most $m-(t+2)-(a-1)$. We are done as
  $a\ge 1$.

  It remains to deal with the case where $a=0$. Recall that $t\ge 3$ and $t=2a+3$, $t=2a+4$, or
  $t=a+2$. This forces $t=3$ or $t=4$. Consequently, $m=6k+3$ when $t=3$, or
  $t=4$ and $m=8k+4$.

  Since $\bv_{r+1}\in \Lambda(\bu^{*k})$, there exist $c_i$'s in $\F_{2^m}$
  such that
  \[\bv_{r+1} = \sum_{j=1}^{t-1}\sum_{i=0}^{k-1} c_{j+ti} \bv_{j+ti}.\]
  It follows that
  \[
  \bv_{r+2}
  = \sum_{j=1}^{t-1}\sum_{i=0}^{k-1} c_{j+ti}^2 \bv_{j+1+ti}
  = \sum_{i=0}^{k-1}c_{(t-1)+ti}^2\bv_{t(i+1)}+
  \sum_{j=1}^{t-2}\sum_{i=0}^{k-1} c_{j+ti}^2 \bv_{j+1+ti}.
  \]
  Since $t>2$ we have $\bv_{r+2}\in\Lambda(\bu^{*k})$.  Note that also
  $\sum_{j=1}^{t-2}\sum_{i=0}^{k-1} c_{j+ti}^{2} \bv_{j+1+ti}\in
  \Lambda(\bu^{*k})$ and thus,
  $\sum_{i=0}^{k-1}c_{t-1+ti}^2\bv_{t(i+1)}\in\Lambda(\bu^{*k})$. However,
  we also have
  $\sum_{i=0}^{k-1}c_{t-1+ti}^2\bv_{t(i+1)}\in\Lambda(\bzero*\bbu^{*k})$.
  Now observe that $V$ is spanned by the $r$ linearly independent vectors
  $\bv_1, \bv_2, \ldots, \bv_{kt}$, $\Lambda(\bzero*\bbu^{*k})$ is spanned
  by the $k$ vectors $\bv_t, \bv_{2t}, \ldots, \bv_{kt}$, and
  $\Lambda(\bu^{*k})$ is spanned by the $r-k$ vectors in $\{ \bv_1,
  \bv_2, \ldots, \bv_{kt}\}\setminus\{ \bv_t, \bv_{2t}, \ldots, \bv_{kt}\}$.
  Therefore, $\sum_{i=0}^{k-1}c_{t-1+ti}^2\bv_{t(i+1)} \in
  \Lambda(\bu^{*k})\cap\Lambda(\bzero*\bbu^{*k})$ has to be the zero vector
  in $\F_{2^m}^r$. This forces $c_{t-1}=c_{(t-1)+t}=\cdots
  =c_{(t-1)+(k-1)t}=0$.

  If $t=4$, then by applying a similar argument on $\bv_{r+3}$, we see that
  $c_{t-2}=c_{(t-2)+t}=\cdots =c_{(t-2)+(k-1)t}=0$. Thus, in both cases, we
  obtain
  \[\bv_{r+1} = \sum_{i=0}^{k-1} c_{1+ti} \bv_{1+ti}.\]

  Let $h$ be the largest integer such that $c_{1+th}\neq 0$. If $h\neq k-1$,
  then $\bv_{r+1+(k-h-1)t}=\sum_{i=0}^{k-1} c_{1+ti}'\bv_{1+ti}$ where
  $c_{1+t(k-1)}'=c_{1+th}^{2^{t(k-h-1)}}\neq 0$. Hence, it follows that
  \[
  V = \begin{cases}
    \Lambda(\bu^{*(k-1)} * (0,0,1) * \bbu^{*(k-h-1)} * (1,1,0) * \bbu^{*h})
    & \text{if $t=3$}\\
    \Lambda(\bu^{*(k-1)} * (0,0,1,1) * \bbu^{*(k-h-1)} * (1,1,0,0) *
    \bbu^{*h}) & \text{if $t=4$}.
  \end{cases}
  \]
  If $h=k-1$, i.e., $c_{1+t(k-1)}\neq 0$, then we see that
  \[
  V = \begin{cases}
    \Lambda(\bu^{*(k-1)} * (0,1,0) * (1,0,1)*\bbu^{*(k-1)})
    & \text{if $t=3$}\\
    \Lambda(\bu^{*(k-1)} * (0,1,0,1) *(1,0,1,0)* \bbu^{*(k-1)} )
    & \text{if $t=4$}.
  \end{cases}
  \]
  When $k\geq 2$, after dropping the zero in the first copy of $\bu$ and the
  last zero in the last copy of $\bbu$, we obtain a vector we require in
  each case. When $k=1$, we deduce that $t=4$ as $m=9$ is excluded
  (cf.~Example~\ref{example:m=9}), and clearly the required vector is then
  $(1,0,1,1,0,1)$.
  \vspace{0.1in}

  {\bfseries Case \boldmath$\Lambda(\ba * \bu^{*k} * \bu^{*b})= V\text{ for
      some } b\le k$:} We consider two subcases depending on the increase in
  dimension in the nested sequence $\Lambda(\ba*\bu^{*k})\subsetneq
  \Lambda(\ba*\bu^{*(k+1)})\subsetneq \cdots \subsetneq V$.

  {\slshape Subcase \boldmath$b<k$:} The dimension of one of the subspaces in
  the sequence increases by more than one compared to that of its predecessor,
  and $1\le b\le k-1$. Thus, $\ba * \bu^{*(k+b)}$ is a vector of length $a+(k+b)t\le
  a+(2k-1)t=m-a-(t+3)$. By construction, this vector does not have more than
  $t-1$ consecutive 1's and we are done.

  {\slshape Subcase \boldmath$b=k$:} The dimension of each vector space in
  the nested sequence $\Lambda(\ba*\bu^{*k})\subsetneq \cdots \subsetneq
  \Lambda(\ba*\bu^{*(2k)})=V$ increases by exactly one compared
  to that of its predecessor. Hence there is a smallest index $j$, $1\le j\le t-1$,
  such that $\bv_{a+kt+j}\notin \Lambda(\ba *\bu^{*k})$ and
  $\Lambda(\ba*\bu^{*(k+1)})= \Lambda(\ba *\bu^{*k})+\F_{2^m}
  \bv_{a+kt+j}$. It follows from
  Lemma~\ref{lemma:nestedsequence}, part~\ref{lemma:nestedsequence:ii},, that
  $V=\Lambda(\ba*\bu^{*2k})= \Lambda(\ba *\bu^{*(2k-1)})+\F_{2^m}
  \bv_{a+(2k-1)t+j}$. Therefore, we conclude
  \[ V=\Lambda(\ba*\bu^{*(2k-1)}*(\underbrace{0,\ldots,0}_{j},1)).\]
  Note that the length of the vector
  $\ba*\bu^{*(2k-1)}*(\overbrace{0,\ldots,0}^{j},1)$ is $r+(k-1)t+(j+1)$,
  which is at most $m-(t+2)-(a-j)$. By construction, it does not have more
  than $t-1$ consecutive 1's. Therefore, we are done if $j\le a$.

  We still have to deal with the case where $j>a$. In this case,
  \[ \Lambda(\ba*\bu^{*k})=\Lambda(\ba*\bu^{*k} *
  (0,\underbrace{1,\ldots,1}_{j-1}))\]
  since $j$ was the smallest index such that $\Lambda(\ba
  *\bu^{*(k+1)})=\Lambda(\ba*\bu^{*k})+ \F_{2^m} \bv_{r+j}$.
  However, it is clear that the set $ \{ \bv_j, \ldots, \bv_{r+(j-1)}\}
  \setminus \{ \bv_{a+t}, \ldots, \bv_{a+kt}\}$ is linearly independent.
  Therefore,
  \[ \Lambda(\ba*\bu^{*k}) = \Lambda((\underbrace{0,\ldots,0}_{j},
  \underbrace{1,\ldots,1}_{t+a-j})*\bu^{*(k-1)} *
  (0,\underbrace{1,\ldots,1}_{j-1}))\]
  as both spaces have the same dimension. It follows that
  \[ V=\Lambda(\ba*\bu^{*(2k-1)}*(0,\underbrace{1,\ldots,1}_j))=\Lambda((\underbrace{0,\ldots,0}_{j},\underbrace{1,\ldots,1}_{t+a-j})*\bu^{*(2k-2)}*(0,\underbrace{1,\ldots,1}_j)).\]
  Deleting leading and tailing zeros, we obtain a vector that has length at
  most $m-(t+2)-a$. This completes our proof.
\end{proof}

Combining Theorem~\ref{theorem:p1} with known constructions, we
have
\begin{theorem}
  \label{theorem:p:1}
Let $m\geq 5$ but $m\neq 9$. Then the largest $d$ of a
non-Denniston maximal arc of degree $2^d$ in $\PG(2,2^m)$
generated by a $\{p,1\}$-map via Theorem~\ref{theorem:Mathon:pq}
is $\floor{\frac{m}{2}}+1$.
\end{theorem}

\begin{proof} Let $p(x)=\sum_{i=0}^{d-1} a_i x^{2^i-1}\in\F_{2^m}[x]$.
Assume that $\trace(p(\lambda))=1$ for all $\lambda\in
A\setminus\{ 0\}$, where $A$ is an additive subgroup of
$\F_{2^d}$. If $m\geq 5$ but $m\neq 9$, and if $m>d>\frac
{m}{2}+1$, then by Theorem~\ref{theorem:p1}, $p(x)$ is a linear
polynomial, hence the maximal arc generated by the $\{p,1\}$-map
is a Denniston maximal arc. This shows that when $m\geq 5$ but
$m\neq 9$, the largest $d$ of a non-Denniston maximal arc of
degree $d$ in $\PG(2,2^m)$ generated by a $\{p,1\}$-map via
Theorem~\ref{theorem:Mathon:pq} is $\leq \floor{\frac{m}{2}}+1$.

On the other hand, there are always $\{p,1\}$-maps generating
non-Denniston maximal arcs of degree $2^{\floor{\frac{m}{2}}+1}$
if $m\geq 5$ (see \cite{Mat02}, \cite{HamMat}, \cite{FLX}). The
conclusion of the theorem now follows.
\end{proof}

We remark that when $m=9$, there is an example of $\{p,1\}$-maps
that generates a non-Denniston maximal arc of degree $2^6$. This
example appears in \cite{HamMat}.

\begin{example}[\cite{HamMat}]
\label{example:m=9} Let $g$ be a primitive element in $\F_{2^9}$.
Note that $73\cdot (2^3-1)=2^9-1$, so $b=g^{73}$ is a primitive
element in $\F_{2^3}$. Let $\mu_i=b^i$ for $i=0,1,2$ and $A=\{
x\in\F_{2^9}\mid \trace(\mu_i x)=0,\ \forall i=0,1,2\}$. That is,
$A=\{x\in\F_{2^9}\mid \trace_{2^9/2^3}(x)=0\}$ since
$\mu_1,\mu_2,\mu_3$ are linearly independent over $\F_2$. Let
$p(x)=x^7 + 1$. Direct computations show that
$\trace(p(\lambda))=1$ for all $\lambda\in A\setminus\{0\}$.
Therefore the set of points on the conics in
$\{F_{p(\lambda),1,\lambda}\mid \lambda\in A\setminus\{0\}\}$
together with the common nucleus $F_0$ forms a non-Denniston
maximal arc of degree $2^6$.
\end{example}

\section{Upper Bound for the Degree of non-Denniston Maximal Arcs in $\PG(2,2^m)$ Generated by $\{p,q\}$-Maps}
\label{section:p-q}

In this section we try to extend the result in previous section to
$\{p,q\}$-maps, where $q$ is not necessarily $1$.

\begin{theorem}
  \label{theorem:pq}
  Let $A$ be an additive subgroup of size $2^d$ in $\F_{2^m}$, and let
  $p(x)=\sum_{i=0}^{d-1} a_ix^{2^i-1}\in\F_{2^m}[x]$, $q(x)=\sum_{i=0}^{d-1}
  b_ix^{2^i-1}\in\F_{2^m}[x]$. Assume that
  $m\geq 7$ but $m\neq 9$, and $m>d>\frac {m} {2}+2$. If
  $\trace(p(\lambda)q(\lambda))=1$ for all $\lambda\in A\setminus\{ 0\}$,
  then $a_2=a_3=\cdots =a_{d-1}=0$ and $b_2=b_3=\cdots =b_{d-1}=0$.
  That is, $p(x)$ and $q(x)$ are both linear and the maximal arc obtained
  from the $\{p,q\}$-map via Theorem~1.2 is a Denniston maximal arc.
\end{theorem}
\begin{proof} Let $r=m-d$. As in the proof of Theorem~\ref{theorem:p1},
we assume that the defining equation of $A$ is
  \[
  \prod_{i=1}^r (1+\trace(\mu_i x))=1,
  \]
where $\mu_i\in\F_{2^m}^*$ are linearly independent over $\F_{2}$.
Also as argued in the proof of Theorem~\ref{theorem:p1}, we may
assume that $\trace(a_0 b_0)=1$. Then
  \begin{align}
    \label{equation:traceproduct:pq}
    \trace(p(x)q(x) +a_0 b_0) \prod_{i=1}^{r}
    (1+\trace(\mu_i x)) & \equiv 0\pmod{x^{2^m}-x}
  \end{align}
  For convenience, set $T(x)=\trace(p(x)q(x)+a_0b_0)$ and
  $S(x)=\prod_{i=1}^{r}(1+\trace(\mu_i x))$. Also as before
  denote the coefficient
  of $x^{2^{i_1}+2^{i_2}+\cdots + 2^{i_s}}$ in $S(x)$ by
  $c(i_1,i_2,\ldots ,i_s)$, where $1\le s\le r$ and
  $m-1\ge i_1>i_2>\cdots >i_s\ge 0$. The remarks about
  $c(i_1,i_2,\ldots ,i_s)$ in the course of proving
  Theorem~\ref{theorem:p1} are of course valid here.

  In the proof of Theorem~\ref{theorem:p1} we use the
  fact that the exponent of any term in the expansion of $\trace(p(x))$
  is a cyclic shift of $(2^i-1)$ for some $i$. This is no longer true for
  $T(x)=\trace(p(x)q(x)+a_0b_0)$ if $q(x)$ is not a constant.
  Instead, the exponent of
  any term in the expansion of $T(x)$ is $2^s((2^j-1)+(2^k-1))$,
  where $m-1\geq s\geq 0$, $d-1\ge j\ge k\ge 0$.
  If $k\ge 1$ then the binary representation of $2^s((2^j-1)+(2^k-1))$
  is a cyclic shift of
  \[
  0\ldots 01\overbrace{0\ldots 0}^{j-k}\overbrace{1\ldots 1}^{k-1}0.
  \]
  The number of 1's in this representation is $k$. If $k=0$ then it is a
  cyclic shift of
  \[
  0\ldots 0\overbrace{1\ldots 1}^j.
  \]
  The number of 1's is $j$. This shows that the maximum number of 1's
  in the binary representations of such exponents is $d-1$. Note that
  if $k=0$ or $k=j$ then such an exponent is  $2^s(2^i-1)$ for some $s$
  and $i$, hence its binary representation is a shift of $i$ consecutive 1's.

  We want to use techniques similar to those in the proofs of
  Theorem~\ref{theorem:p1}. That is, we will be
  looking at the coefficients of various terms in $T(x)\cdot S(x)$.
  We will be particularly interested in terms $x^e$ in $T(x)$, where the
  exponent $e$ has $(d-1)$ or $(d-2)$ ones in its binary representation.
  If $e$ has $(d-1)$ ones, it must be a shift of $2^d-2$ or $2^{d-1}-1$.  The
  coefficients of $x^{2^d-2}$ and $x^{2^{d-1}-1}$ in $T(x)$ are
  $a_{d-1} b_{d-1}$ and $a_0 b_{d-1}+ a_{d-1} b_0$,
  respectively. If $e$ has $(d-2)$ ones, then it must be a shift of one
  of $2^{d-2}-1$, $2^{d-1}-2$, $2^{d-1}+2^{d-2}-2$. The coefficients of
  $x^{2^{d-2}-1}$, $x^{2^{d-1}-2}$, $x^{2^{d-1}+2^{d-2}-2}$ are
  $a_0b_{d-2} + a_{d-2} b_0$, $a_{d-2} b_{d-2}$ and $a_{d-1} b_{d-2} + a_{d-2} b_{d-1}$, respectively.

  {\bfseries Claim: \boldmath$a_{d-2} b_{d-1}+a_{d-1}b_{d-2}=0$}. Consider
  the coefficient of $x^{2^m-2^{d-2}-2}$ in $T(x)\cdot S(x)$. The binary
  representation of the exponent is
  \[
  \overbrace{1\ldots 1}^{r}10\overbrace{1\ldots 1}^{d-3}0,
  \]
  which has $(m-2)$ ones. The maximum number of 1's in the exponent of
  any summand in $S(x)$ is $r$ and the maximum number of 1's in the
  exponent of  any summand in $T(x)$ is $d-1$. When adding two exponents
  (written in their binary representations), any carry that may occur reduces
  the number of 1's in the sum. Since we are interested in an
  exponent whose number of 1's is $(m-2)$, it can only be obtained as a sum
  of two exponents (one is the exponent of a summand in $T(x)$, the other
  in $S(x)$) with at most one carry.

  Suppose the exponent $2^m-2^{d-2}-2$ is obtained without carry.
  Using the assumption that $d>\frac {m}{2}+2$, we have $d-3>r+1$.
  So there is only one possibility.
  \begin{align*}
    \overbrace{1\ldots 1}^{r}10\overbrace{1\ldots 1}^{d-3}0 & =
    \overbrace{1\ldots 1}^{r}000\ldots 00 + 0\ldots 010\overbrace{1\ldots1}^{d-3}0.
  \end{align*}
  Hence  $0\ldots 010\overbrace{1\ldots1}^{d-3}0$ must come from the
  exponent of $x^{2^{d-1}+2^{d-2}-2}$ in $T(x)$, whose coefficient is
  $a_{d-2}b_{d-1}+a_{d-1}b_{d-2}$, and
  $\overbrace{1\ldots 1}^{r}000\ldots 00$ must come from
  $x^{2^{m-1}+\cdots +2^{d}}$ in $S(x)$, whose coefficient is
  $c(m-1,m-2,\ldots, d)$.

  Now suppose that the exponent $2^m-2^{d-2}-2$ is obtained with a carry,
  which means that the contribution from $T(x)$ is a
  shift of $2^{d-1}-1$. Then it has to be exactly one carry which has to
  occur at position $d-2$ since $d-3>r+1$. There is no way of realizing
  this with any shift of $2^{d-1}-1$.

  Therefore the coefficient of $x^{2^m-2^{d-2}-2}$ in $T(x)\cdot S(x)$ is
  $(a_{d-2} b_{d-1}+a_{d-1}b_{d-2})\cdot c(m-1, m-2, \ldots, d)$, and by
  (\ref{equation:traceproduct:pq}), we have
  \begin{align*}
    (a_{d-2} b_{d-1}+a_{d-1}b_{d-2})\cdot c(m-1, m-2, \ldots, d)
    & = 0.
  \end{align*}
  Noting that $c(m-1,m-2,\ldots, d)$ is a Moore determinant, which is nonzero,
  we conclude that $a_{d-2} b_{d-1}+a_{d-1}b_{d-2}=0$.

  After proving the above claim, observe that now the exponent of any term in
  $T(x)$ whose number of 1's is $d-1$ or $d-2$ has to
  be a cyclic shift of $2^{d-1}-1$ or $2^{d-2}-1$. Thus, we are
  ready to proceed as in the proof of Theorem~\ref{theorem:p1}.

  {\bfseries Claim:
    \boldmath$a_0^2 b_{d-2}^2+a_{d-2}b_{d-2}+b_{0}^2a_{d-2}^2=a_0^2
    b_{d-1}^2+a_{d-1}b_{d-1}+b_{0}^2a_{d-1}^2=0$}. The coefficient of
  $x^{2^{d-1}-2}$ in $T(x)$ is
  \begin{align}
    \label{equation:coefficient:pq:1}
    & a_0^2 b_{d-2}^2+a_{d-2}b_{d-2}+b_{0}^2a_{d-2}^2
    \intertext{and the coefficient of $x^{2(2^{d-1}-1)}$ is}
    \label{equation:coefficient:pq:2}
    & a_0^2 b_{d-1}^2+a_{d-1}b_{d-1}+b_{0}^2a_{d-1}^2.
  \end{align}
  Considering the coefficient of $x^{2^m-2^{d-1}-2}$ and that of
  $x^{2^m-2^d-2}$ in $T(x)\cdot S(x)$, we
  obtain equations similar to (\ref{equation:coefficient:1}) and
  (\ref{equation:coefficient:2}) with the expressions in
  (\ref{equation:coefficient:pq:1}) and (\ref{equation:coefficient:pq:2})
  taking the place of $a_{d-2}$ and $a_{d-1}$ in
  (\ref{equation:coefficient:1}) and (\ref{equation:coefficient:2})
  respectively. Thus, using the same reasoning as in the proof of
  Theorem~\ref{theorem:p1}, our claim follows.

  {\bfseries Claim: \boldmath$a_{d-1}=a_{d-2}=b_{d-1}=b_{d-2}=0$}. Since
  $\trace(a_0 b_0)=1$, the binary quadratic form $a_0^2x^2+xy+b_0^2y^2$ over
  $\F_{2^m}$ has only trivial zeros. Therefore, from
  \begin{align*}
    a_0^2 b_{d-2}^2+a_{d-2}b_{d-2}+b_{0}^2a_{d-2}^2 & = 0\\
    a_0^2 b_{d-1}^2+a_{d-1}b_{d-1}+b_{0}^2a_{d-1}^2 & = 0
  \end{align*}
  we obtain $a_{d-2}=b_{d-2}=0$ and $a_{d-1}=b_{d-1}=0$.

  {\bfseries Claim: \boldmath$a_{d-3}=\cdots =a_{r+1}=b_{d-3}=\cdots
    =b_{r+1}=0$}. Let $d-2>k>r$ and suppose that $a_j=b_j=0$ for $j>k$.
  Consider the coefficient of $x^{2^m-2^d+2^{k+1}-2}$ in $T(x)\cdot S(x)$.
  The exponent of this monomial has
  binary representation
  \[
  \overbrace{1\ldots 1}^{r}\overbrace{0\ldots 0}^{d-k-1}\overbrace{1\ldots 1}^k 0,
  \]
  which has $(m-1)$ ones. This exponent can only be obtained as a sum
  of two exponents (one is the exponent of a summand in $T(x)$, the other
  in $S(x)$) without carry. As we discussed previously, there are three
  ways such that the number of 1's in the binary representation
  of $2^j-1+2^k-1$ is $k>0$. These are $2^k-1+2^0-1$ (the coefficient of
  $x^{2^k-1+2^0-1}$ in $T(x)$ is $a_k b_0+a_0 b_k$), $2^k-1+2^k-1$
  (the coefficient of $x^{2^k-1+2^k-1}$ in $T(x)$ is $a_k b_k$), and
  $2^j-1+2^k-1$ where $j>k$. In the last case, the coefficient of
  $x^{2^j-1+2^k-1}$ is $\sum_{j>k}(a_k b_j + b_k a_j)$, which is zero
  since $a_j=b_j=0$ for $j>k$.

  Hence the coefficient of $x^{2^m-2^d+2^{k+1}-2}$ in $T(x)\cdot S(x)$ is
  \begin{align*}
    (b_0^2 a_k^2 + a_k b_k + a_0^2 b_k^2)\cdot c(m-1, m-2, \ldots, d) & = 0.
  \end{align*}
  As before, $c(m-1, m-2, \ldots, d)$ is a Moore determinant, which is nonzero.
  Therefore $(b_0^2 a_k^2 + a_k b_k + a_0^2 b_k^2)=0$. Since
  $\trace(a_0b_0)=1$, we have $a_k=b_k=0$.

  Note that in the case where $d=m-1$, the above claims already
  show that $a_2=a_3=\cdots =a_{d-1}=0$ and $b_2=b_3=\cdots
  =b_{d-1}=0$, so $p(x)$ and $q(x)$ are both linear. Also observe
  that when $m=7$ (resp. 8), the only admissible $d$ is 6 (resp.
  7). In both cases, $m-d=1$, so $p(x)$ and $q(x)$ are both
  linear. Hence from now on, we will assume that $m\geq 10$ and
  $m-1>d>\frac {m}{2}+2$.

  {\bfseries Claim: \boldmath$a_{r}=\cdots =a_3=b_{r}=\cdots =b_3=0$}.
  Let $3\le t\le r$ and assume that $a_j=b_j=0$ for $j>t$. Since
  $r\le \frac{m-3}{2}$ and $m\geq 10$, by Theorem~\ref{lemma:gap},
  there exist $0=i_1<i_2<\cdots < i_{r}\le m-t-3$ such that $c(i_1,i_2,\ldots
  ,i_{r})\neq 0$ and the number of consecutive 1's in $\{ i_1,i_2,\ldots
  ,i_{r}\}$ is at most $t-1$.  Now we consider the exponent $2^m-2^{m-t} +
  \sum_{j=1}^{r}2^{i_j}$ and we see that it can only be obtained in one way
  as a sum of two exponents, one from $T(x)$, the other from $S(x)$.
  \begin{align*}
    0\overbrace{1\ldots 1}^{k}0\overbrace{\ldots 1\ldots 1}^{m-k-2} & =
    0\overbrace{0\ldots 0}^{k}0\overbrace{\ldots 1\ldots 1}^{m-k-2} +
    0\overbrace{1\ldots 1}^{k}0\overbrace{0\ldots 0}^{m-k-2}.
  \end{align*}
  It follows from (\ref{equation:traceproduct:pq}) that
  \begin{align*}
    (b_0^2 a_k^2 + a_k b_k + a_0^2 b_k^2)\cdot c(i_1,i_2,\ldots, i_{r}) &
    = 0
  \end{align*}
  and hence, $a_k=b_k=0$.

  {\bfseries Claim: \boldmath$a_2=b_2=0$}. As in
  Theorem~\ref{theorem:p1} we consider the
  quadratic form $Q(x)=\trace(p(x)q(x)+a_0 b_0)$ over $V=\F_{2^m}$. Note
  that since $\trace(a_0b_0)=1$, the assumption that
  $\trace(p(\lambda)q(\lambda))=1$ for all $\lambda\in A\setminus\{ 0\}$
  implies that $Q(\lambda)=0$ for all $\lambda\in A$, where $\abs{A}=2^d$. The
  bilinear form associated with $Q(x)$ is
  \begin{align*}
    B(x,y) & = \trace\left((a_0^2b_2^2+a_2 b_2+a_2^2b_0^2)(x y^2+y x^2)^2\right).
  \end{align*}
  \begin{align*}
    \rad V & = \{ x\in V\mid \trace\left((a_0^2b_2^2 + a_2 b_2 + a_2^2b_0^2) (x y^2+yx^2)^2\right) = 0,\,\forall_{y\in V}\}\\
    & = \{ x\in V\mid x^3 = (a_0^2 b_2^2 + a_2 b_2 + a_2^2b_0^2)^{-\frac{1}{2}}\}\cup \{0\}.
  \end{align*}
  As discussed in the proof of Theorem~\ref{theorem:p1}, if
  $a_0^2 b_2^2 + a_2 b_2 + a_2^2b_0^2\neq 0$, then the maximum dimension of
  a subspace of $V$ on which $Q$ vanishes is at most
  $\floor{\frac{m}{2}}+1$. But we knew that $Q(x)$ vanishes on $A$, which
  has $\F_2$-dimension $d$, and $d>\frac{m}{2}+2$. This is a contradiction.
  Hence $a_0^2 b_2^2 + a_2 b_2 + a_2^2 b_0^2=0$. Combining this with $\trace(a_0b_0)=1$,
  we obtain $a_2=b_2=0$.

  So we have proven that both $p(x)$ and $q(x)$ must be linear,
  by the last part of Theorem~\ref{theorem:Mathon:pq}, the maximal arc
  generated by this $\{p,q\}$-map is a Denniston maximal arc. This completes
  the proof.
\end{proof}

Combining Theorem~\ref{theorem:pq} with known constructions in
\cite{Mat02}, \cite{HamMat} and \cite{FLX}, we have
\begin{theorem}
  Let $m\geq 7$ but $m\neq 9$. Then the largest $d$ of a non-Denniston
  maximal arc of degree $2^d$ in $\PG(2,2^m)$ generated by a $\{p,q\}$-map
  via Theorem~\ref{theorem:Mathon:pq} is either $\floor{\frac{m}{2}}+1$ or
  $\floor{\frac{m}{2}}+2$.
\end{theorem}

It is an interesting question whether there exists a $\{p,q\}$-map
generating a non-Denniston maximal arc in $\PG(2,2^m)$ of degree
$\floor{\frac{m}{2}}+2$ when $m\geq 7$. We remark that in the case $m=5$,
there is an example of $\{p,q\}$-maps which generates a non-Denniston
maximal arc of degree $16$ in $\PG(2,32)$ (\cite[p.~362]{Mat02}).

\vspace{0.1in}

\noindent{\bf Acknowledgements:} We thank Henk D. L. Hollmann for
useful discussions on Mathon's construction of maximal arcs. Part
of this work was carried out during a visit of the third author to
National University of Singapore. The third author would like to
thank Department of Mathematics, National University of Singapore
for its hospitality. The research of the third author was also
partially supported by NSA grant MDA 904-03-1-0095.

\end{document}